%%%%%%%%%%%%%%%%%%%%%%%%%%%%%%%%%%%%%%%%%%%%%%%%%%%%%%%%%%%%%%%%%%%%%%%%%%%
%% Martino, John; Priddy, Stewart
%% 
%% A Classification of the Stable Type of $BG$
%% 
%% We give a classification of the $p$--local stable homotopy type of $BG$,
%%   where $G$ is a finite group, in purely algebraic terms. $BG$ is 
%%   determined by conjugacy classes of homomorphisms from $p$--groups into
%%   $G$. This classification greatly simplifies if $G$ has a normal Sylow 
%%   $p$--subgroup; the stable homotopy types then depends only on the Weyl
%%   group of the Sylow $p$--subgroup. If $G$ is cyclic mod $p$ then $BG$ 
%%   determines $G$ up to isomorphism. The last class of groups is 
%%   important because in an appropriate Grothendieck group $BG$ can be 
%%   written as a unique linear combination of $BH$'s, where $H$ is cyclic 
%%   mod $p$.
%% 
%% publ:  Bull. Amer. Math. Soc. (N.S.) 27(1992) no. 1
%% pp:    165-170
%% type:  Research Announcement        markup: amstex    file size: 20K
%% 
%% copyright: American Math. Society copyright; see end of article
%% 
%% Include files necessary for this article: bull-ppt.tex
%% 
%%%%%%%%%%%%%%%%%%%%%%%%%%%%%%%%%%%%%%%%%%%%%%%%%%%%%%%%%%%%%%%%%%%%%%%%%%%
\input amstex
\documentstyle{amsppt}
\input bull-ppt
\keyedby{bull300e/pah}
%\pagewidth{6.0 truein}
%\pageheight{8.7 truein}

\define\lamar{\mathop{\longrightarrow}}
\define\lams{\mathop{\rightarrow}}
\topmatter
\cvol{27}
\cvolyear{1992}
\cmonth{July}
\cyear{1992}
\cvolno{1}
\cpgs{165-170}
%\ratitle
\title A Classification of the Stable Type of $BG$\endtitle
\author John Martino and Stewart Priddy\endauthor
\thanks The first author was partially supported by NSF 
Grant 
DMS-9007361, the second by NSF Grant DMS-880067 and the 
Alexander von
Humboldt Foundation\endthanks
\date November 20, 1991\enddate
\subjclass Primary 55R35; Secondary 20J06, 
55P42\endsubjclass
\address Department of Mathematics, University of Virginia, 
Charlottesville, Virginia 22903\endaddress
\address Department of Mathematics, Northwestern 
University, Evanston, Illinois 60208\endaddress 
\abstract We give a classification of the $p$--local 
stable homotopy type of 
$BG$, where $G$ is a finite group, in purely algebraic 
terms.  $BG$ is 
determined by conjugacy classes of homomorphisms from 
$p$--groups into $G$.
This classification greatly simplifies if $G$ has a normal 
Sylow 
$p$--subgroup; the stable homotopy types then depends only 
on the Weyl group 
of the Sylow $p$--subgroup.  If $G$ is cyclic mod $p$ then 
$BG$ determines $G$ 
up to isomorphism.  The last class of groups is important 
because in an 
appropriate Grothendieck group $BG$ can be written as a 
unique linear 
combination of $BH$'s, where $H$ is cyclic mod 
$p$.\endabstract
\endtopmatter
 
\heading 0. Introduction and statement of main 
results\endheading

Let $G$ be a finite group. In this note we give a 
classification of the stable homotopy type of $BG$ in 
terms of $G$. Our 
analysis shows that for each prime number $p$, the 
$p$-local stable type of 
$BG$ depends on the homomorphisms from $p$-groups $Q$ into
$G$.

The suspension spectrum of $BG$ and, in particular, its 
wedge summands have 
played an important role in homotopy theory. In a previous 
paper \cite{MP}, 
the authors have given a characterization of the 
indecomposable summands of 
$BG$ in terms of the modular representation theory of 
$\roman{Out}(Q)$ 
modules for 
$Q<P$ the Sylow subgroup of $G$. It is this 
characterization that we use to 
study the stable type of $BG$. For another such 
characterization see 
\cite{BF}. 
 
It is known that the stable type of $BG$ does not 
determine $G$ up to 
isomorphism. A simple example (due to  N. Minami) is given 
by $Q_{4p} \times 
\bold Z/2$ and $D_{2p} \times \bold Z/4$  where $p$ is an 
odd prime, $Q_{4p}$
is the generalized quaternion group \cite{CE} of order 
$4p$, and $D_{2p}$ is 
the dihederal group of order $2p$. The situation is even 
worse for  $p$-local
classifying spaces since $BG$ and $BG/O_{p'}(G)$ have 
isomorphic mod $p$ 
homology and hence equivalent stable types. Here 
$O_{p'}(G)$ is the maximal 
normal subgroup of $G$ of order prime to $p$. However 
there is a positive 
result in this direction, due to  Nishida \cite{N}, who 
established the 
following: Suppose $G_{1}$, $G_{2}$ are  finite groups 
with Sylow
$p$-subgroups
$P_1, P_2$, then $BG_1 \simeq BG_2$ stably at $p$ implies 
$P_1 
\approx P_2$. Our main result is a necessary and sufficient 
condition. 

\proclaim{Theorem 0.1 {\rm (Classification)}} For two 
finite groups
$G_1, G_2$ the following are equivalent\RM:

{\rm (1)} Localized at $p$, $BG_1$ and $BG_2$ are stably 
homotopy equivalent.

{\rm (2)} For every $p$-group $Q$,
$$\bold F_p \roman{Rep}(Q,G_1) \approx \bold F_p 
\roman{Rep}(Q,G_2)$$
as $\roman{Out}(Q)$ modules.  $\roman{Rep}(Q,G) = 
\roman{Hom}(Q,G)/G$ with $G$ acting by conjugation.

{\rm (3)} For every $p$-group $Q$,
$$\bold F_p \roman{Inj}(Q,G_1) \approx \bold F_p 
\roman{Inj}(Q,G_2)$$
as $\roman{Out}(Q)$ modules.  $\roman{Inj}(Q,G) < 
\roman{Rep}(Q,G)$ consists of conjugacy classes of 
injective homomorphisms.
\endproclaim 

Nishida's theorem follows since the largest $Q$ for which 
$\roman{Inj}(Q,G)$ is 
nonzero is the Sylow $p$-subgroup of $G$.  We will refer 
to the common Sylow
$p$-subgroup as $P$. It should not be concluded from (3) 
that the $G_1$ 
conjugacy classes of a $p$-subgroup correspond to the 
$G_2$ conjugacy classes, 
although the number of classes is equal.

An important application of Theorem 0.1 is to the case of 
a normal Sylow
$p$-subgroup, in particular, the cyclic mod $p$ groups.

\dfn{Definition 0.2} Two subgroups $H,K<G$ are called \it 
pointwise 
conjugate in G\ \rm if there is a bijection of sets $H
\lamar\limits^\alpha K$ such that $\alpha (h) = 
g_{h}^{-1}h g_h$  for 
$g_h\in G$ depending on $h\in H$. \enddfn

Alternately it is easy to see that an equivalent condition 
is 
$$
               |H \cap(g)| = |K \cap(g)|
$$
for all $g\in G$, where $(g)$ denotes the conjugacy class 
of $g$.

Let $W_G(H)$ denote the \it{Weyl group}\rm\ of $H<G$, i.e.,

$$
                    W_G(H)=N_G(H)/H\cdot C_G(H)
$$
where $N_G(H)$  is the normalizer and $C_G(H)$  is the 
centralizer of $H$ in 
$G$. Then $W_G(H)<\roman{Out}(H)$. 

\proclaim{Theorem 0.3} Suppose $G_1, G_2$ are finite 
groups with normal Sylow
$p$-subgroups $P_1, P_2$. Then $BG_1$ and $BG_2$ have the 
same stable
homotopy
type, localized at $p$, if and only if $P_1 \approx P_2$ 
$(\approx P$ say\RM)
and $W_{G_1}(P)$ is pointwise conjugate to $W_{G_2}(P)$ in 
$\roman{Out}(P)$. 

\endproclaim

\dfn{Definition 0.4} $G$ is called \it cyclic mod $p$ \rm 
(or
$p$-hypoelementary) if a Sylow $p$-subgroup 
$P$ is normal and $C=G/P$ is a cyclic $p'$-group, i.e., 
has order prime to
$p$.
We say $G$ is \it{reduced}\ \rm if $O_{p'}(G)=1$. 
\enddfn

For a cyclic mod $p$ group $G$, being reduced is 
equivalent to $W_{G}(P)=C$.
 
\proclaim{Theorem 0.5} Suppose $G_1$, $G_2$ are reduced 
cyclic mod $p$
groups.
Then $BG_1 \simeq BG_2$ stably at $p$ if and only if $G_1
\approx G_2$.
\endproclaim

Cyclic mod $p$ groups are important for several reasons. 
Their mod $p$
cohomology is computed as the ring of invariants 
$H^*(G)=H^*(P)^C$. On the 
level of stable homotopy one also has the Minami-Webb 
Formula \cite{M}: Let 
$\scr C_{p}(G)$ be the set of cyclic mod $p$ subgroups of 
$G$. Then 

$$ 
                     BG \simeq \bigvee_{(H)} 
f(H)/[N_G(H):H]\ \ BH 
$$ 
where $H$ runs over the conjugacy classes of $\scr 
C_{p}(G)$ and $f:\scr 
C_{p}(G) \longrightarrow \bold Z$ is the M\"{o}bius 
function given by 
$\sum_{J<K} f(K)=1$
for $J, K\in \scr C_p(G)$.

This expression is to be interpreted in a Grothendieck 
group of spectra,
that is, multiplying through by a common denominator and 
moving the negative 
terms to the left side results in a valid formula in the 
category of
$p$-local 
spectra. 

The next result shows the uniqueness of the Minami-Webb 
Formula.

\proclaim{Theorem 0.6 {\rm (Linear independence of cyclic 
mod $p$ groups)}}
Suppose
$$   
          \bigvee\limits^m_{i=1} \alpha_i BH_i \simeq 
\bigvee\limits^n_{j=1}\alpha'_j BH'_{j} 
$$ 
for $H_i, H_{j}'$ reduced cyclic $\mod p$ groups 
with $\alpha_i,\alpha'_j \in \bold Q$. Then $m=n$ and up 
to a permutation of 
indices $\alpha_i=\alpha'_j, H_i=H'_j$. 
\endproclaim

The note is organized as follows.  In \S1 we discuss the 
classification
(Theorem 0.1).  In \S2 we discuss the normal Sylow
$p$-subgroup case and explain how it 
relates to the general case. Section 3 is devoted to 
cyclic mod $p$ groups.

Preliminary material and background references for this 
note may be found 
in \cite{MP}. All spectra are assumed to be localized at 
$p$ and all 
cohomology groups are taken with simple coefficients in 
$\bold F_{p}$
unless otherwise noted.

\heading 1. Classification theorem\endheading

In this section we briefly discuss Theorem 0.1.

(1) implies (2) since $\bold F_p \roman{Rep}(Q,G)$ is a 
quotient of $\{ BQ,BG \}
\otimes 
\bold F_p$ and the quotient is preserved by stable 
homotopy equivalences.

The equivalence of (2) and (3)
follows from induction on the order of the $p$-group $Q$ 
and the fact that
$$
\bold F_p \roman{Rep}(Q,G)=\bigoplus_R \bold F_p 
\roman{Surj}(Q,R) 
\otimes_{\bold F_p \roman{Out}(R)} \bold F_p 
\roman{Inj}(R,G)$$
where $R$ runs over isomorphism classes of $p$-groups and 
$\roman{Surj}(Q,R)<
\roman{Rep}(Q,R)$ 
consists of conjugacy classes of surjective homomorphisms. 

We show that (3) implies (1) by showing that $BG_1$ and 
$BG_2$ have the same 
indecomposable stable summands with the same 
multiplicities.  If $X$ is such 
an indecomposable summand we show that for certain 
$p$-groups $Q$, $X$ is in 
one-to-one correspondence with a simple module $M$ of a 
quotient ring $R(Q)$
of the outer endomorphism ring of $Q$.  $\bold F_p 
\roman{Out}(Q)$ 
is a subring of 
$R(Q)$. 

To show statement (1) we need another statement equivalent 
to statement (3).
Let $K(Q,G) < \roman{Inj}(Q,G)$ be the classes of all 
injective homomorphisms 
$\alpha$ such that $C_G(\roman{Im}\,\alpha 
)/Z(\roman{Im}\,\alpha )$ 
is a $p'$-group.
An inductive argument shows that the statement 
$$ \bold F_p K(Q,G_1) \approx \bold F_p K(Q,G_2) \tag 4 $$
as $\roman{Out}(Q)$ 
modules, for all $p$-groups $Q$, is equivalent to 
statement (3).

For $\alpha \in\roman{Inj}(Q,G)$, let
$W = N_G(\roman{Im}\,\alpha )/\roman{Im}\,
\alpha$, and $\overline W = \sum_{w \in W} \bar w$,
where $\bar w$ means view $w$ as an element in 
$\roman{Out}(Q)$.
$\overline W \neq 0 \mod p$ if and only if $\alpha \in 
K(Q,G)$.
For conjugacy classes $Q_j$ of $p$-subgroups of $G_1$ (or 
$G_2$), where $X$
is in one-to-one correspondence with a simple module of 
$R(Q_j)$ and the 
inclusion $Q_j \hookrightarrow G_1$ (or $G_2$) is in 
$K(Q_j,G_i)$, $i = 1$ or 
$2$, we find that the multiplicity of $X$ in $BG_1$ (or 
$BG_2$) equals 
$$
\sum_j \dim_k \overline W_j M_j,
$$ 
where $k=\roman{End}(M_j)$.                   
Finally we show that if $W<\roman{Out}(Q)$ then
$$\overline W M = \roman{Hom} (1^{\roman{Out}(Q)}_W,P_M)/
\roman{Hom}(1^{\roman{Out}(Q)}_W,\roman{Ker}\pi(P_M)) $$
where $P_M\lams\limits^\pi M$ is the projective cover of 
$M$ as an 
$\roman{Out}(Q)$ module.  Thus the 
multiplicity depends only on the permutaion modules 
generated by the
$\overline W_j$'s.  $\bold F_p K(Q,G)$ is a direct sum of 
such permutation 
modules, so statement (4) implies statement (1).

\heading 2. Groups with normal Sylow 
$p$-subgroups\endheading

In this section we discuss Theorem 0.3 and its relation to 
 Theorem 0.1. We 
also give an example illustrating Theorem 0.3 and show how 
the stable and 
unstable types of classifying spaces differ.         

In Theorem 0.3 only a condition involving the Sylow 
$p$-subgroup is given 
whereas in Theorem 0.1 all $Q<P$ must be considered. This 
can be explained as 
follows.
Suppose $G$ has a normal Sylow $p$-subgroup $P$. Then 
$G=P\rtimes H$ for a 
$p'$-group $H$ and $BG\simeq B(G/O_{p'}(G)) = B(P\rtimes 
W_{G}(P))$.
Thus we may assume $G=P\rtimes W_{G}(P)$ from which we have 
$\roman{Rep}(Q,G)=\roman{Rep}(Q,P)/W_{G}(P)$, $Q$ a 
$p$-group. Thus 
$$
\bold F_p \roman{Rep}(Q,G)=\bold F_p 
(\roman{Rep}(Q,P)/W_{G}(P))=\bold F_p 
\roman{Rep}(Q,P) \otimes_{\roman{Out}(P)} 1^{\roman{Out}(P)}
_{W_{G}(P)}
$$
where $1^G_H = \bold F_p (G/H)$. Character theory implies 
$W_{G_{1}}(P)$ and $W_{G_{2}}(P)$ are pointwise conjugate 
if and only if the 
induced $\roman{Out}(P)$ modules $1
^{\roman{Out}(P)}_{W_{G_{1}}(P)}$ and 
$1^{\roman{Out}(P)}_{W_{G_{2}}(P)}$ 
are isomorphic. Thus in this case the hypothesis 
of Theorem 0.3 is equivalent to statement (2) of Theorem 
0.1.

\dfn{Example 2.1} Let $p,l$ be different primes and let 
$V$ be an elementary abelian $p$-group of rank $l^{n}$. 
Let $H_{i}$
be two nonisomorphic $l$-groups of exponent $l$ and order 
$l^{n}$ (e.g., for 
$n=3, l>2$, let $H_{1}= (\bold Z/l)^{n}$, $H_{2}= 
U_{3}(\bold F_{l})$, the 
$3\times 3$ upper triangular matrices with $1$'s on the 
diagonal). These 
groups act on $P$ by means of the regular representation 

$$
H_{i}\overset{{\roman{reg}}}\to\hookrightarrow 
\Sigma_{l^n} <
\roman{GL}_{l^n}(\bold F_{p})\.
$$
It follows that 
$G_{i}= V\rtimes H_{i}$
are not isomorphic and satisfy $O_{p'}(G_{i}) = 1$. 
Furthermore $H_{1}$ is 
pointwise conjugate to $H_{2}$ in $\roman{GL}(V) 
= \roman{Out}(V)$. In fact, the nontrivial 
elements of $H_{1}$ and $H_{2}$ are all 
conjugate in $\Sigma _{l^n}$ since each is an $l^{n-1}$ 
fold product of 
disjoint $l$ cycles and the conjugacy of elements is 
determined by their 
cycle structure. By Theorem 0.3, $BG_{1}$ is stably 
equivalent to 
$BG_{2}$ at $p$.

We are indebted to Hans-Werner Henn and Nick Kuhn for 
pointing out that
this example is also interesting because completed at $p$, 
$BG_{1}$ and 
$BG_{2}$ have different (unstable) homotopy types. Thus 
$H^{*}BG_{1}$
and $H^{*}BG_{2}$ are isomorphic in $\scr U$, the category 
of unstable 
modules over the Steenrod algebra $\scr A$, but not in 
$\scr K$, the category
of unstable algebras over $\scr A$. 

To prove this we recall the following result of Adams, 
Miller-Wilkerson,
and Lannes \cite{L}, \cite{HLS}: For $G$ a compact Lie 
group and $V$ an 
elementary abelian $p$-group, the correspondence $f\mapsto 
(Bf)^{*}$
induces a bijection of $\roman{End}(V)$ sets
$$
        \roman{Rep}(V,G)\approx \roman{Hom}_{\scr 
K}(H^{*}BG,H^{*}BV)\.
$$
Thus an isomorphism of unstable algebras $H^{*}BG_{i}$ 
would imply an 
isomorphism of $\roman{End}(V)$ sets $\roman{Rep}(V,G_{i})=
\roman{End}(V)/H_{i}$ and thus
an isomorphism of $\roman{GL}(V)$ sets 
$\roman{GL}(V)/H_{i}$. However, this is false 
since $H_{1}\not\approx H_{2}$.
\enddfn

\heading 3. Cyclic $\mod p$ groups\endheading

In this section we discuss cyclic mod $p$ groups, which 
are determined by
their 
classifying spaces, and that their classifying spaces 
determine
all
other classifying spaces.

Theorem 0.5 is a direct corollary of Theorem 0.3, since 
the generators of the
Weyl groups will be conjugate.

The uniqueness of the expression in the Minami-Webb 
formula follows from a
result in Harris and Kuhn \cite{HK}.

\dfn{Definition 3.1} Let $G$ be a finite group and $p$ a 
prime.

(1) Let $A_p(G)$ be the Burnside ring generated by 
isomorphism classes of 
finite $G$-sets $X$ such that $X^H=\emptyset$ if $H$ is 
not a $p'$-subgroup
of 
$G$, with addition given by disjoint union.

(2) Let $R_p(G)$ be the modular representation ring with 
generators $[M 
\rbrack$, where $M$ is an $\bold{F}_p [G]$-module, and 
relations 
$[M] = [M_1 ]+ [M_2]$ if there is a 
short exact sequence $0 \rightarrow M_1 \rightarrow M 
\rightarrow M_2 
\rightarrow 0$.

(3) Let $\psi_G :A_p(G) \longrightarrow R_p(G)$ be the map 
defined by $\psi_G
(X) = \bold{F}_p [X]$.
\enddfn
\demo{N.B} $A_p(G)$ is generated by $G/H$ where $H$ is a 
$p'$-subgroup 
and
$\psi_G (G/H)=~1_H^G$.
\enddemo

\proclaim{Proposition 3.2 \cite{HK}\RM} For any finite 
group $G$, rank 
$\roman{Im}\,\psi_G$ 
equals the number of conjugacy classes of cyclic 
$p'$-subgroups of
$G$.
\endproclaim

Theorem 0.6 is shown by grading the Grothendieck group of 
classifying spaces 
by Sylow $p$-subgroups, then using the linear independence 
of $p'$-cyclic 
permutaion modules and their correspondence with 
classifying spaces of 
reduced cyclic mod $p$ groups.

\heading Acknowledgments\endheading
We thank David Benson, Nick Kuhn, Guido Mislin, and 
especially Hans-Werner 
Henn for several useful conversations on the topics in 
this note. We are also
grateful to the University of Heidelberg for support 
during part of the time 
this research was done.


\Refs
\ra\key{MMM}
\ref\key BF \by D. Benson and M. Feshbach \paper Stable 
splittings of 
classifying spaces of finite groups \jour Topology 
\toappear \endref 

\ref\key CE \by H. Cartan and S. Eilenberg \book 
Homological algebra \publ
Princeton Univ. Press  \publaddr Princeton, NJ \yr 1956 
\endref 

\ref\key HK \by J. Harris and N. Kuhn \paper Stable 
decompositions of 
classifying spaces of finite abelian $p$-groups \jour 
Math. Proc. Cambridge 
Philos.
Soc. \vol 103 \yr 1988 \pages 427--449 \endref 

\ref\key L \by J. Lannes \paper Sur la cohomologie modulo 
$p$ des
$p$-groupes
\'abeliens \'el\'ementaires \jour Homotopy Theory Proc. of 
the Durham Symp. 
1985, London Math. Soc. Lecture Notes Series, vol. 
117 \publ Camb. Univ. Press\publaddr Cambridge and New 
York, 1987\endref

\ref\key HLS \by H.-W. Henn, J. Lannes, and L. Schwartz 
\paper Analytic
functors,
unstable algebras and cohomology of classifying spaces 
\jour Contemp.
Math., vol. 96,  Amer. Math. Soc.,
Providence, RI, 1989, pp. 197--220 \endref 

\ref\key MP \by J. Martino and S. Priddy \paper The 
complete stable
splitting
for the classifying space of a finite group \jour Topology 
\toappear 
\endref 

\ref\key M \by N. Minami \paper A Hecke algebra for 
cohomotopical Mackey
functors, the stable homotopy of orbit spaces, and a 
M\"{o}bius function
\paperinfo preprint \endref 

\ref\key N \by G. Nishida \paper Stable homotopy type of 
classifying 
spaces of finite groups \jour Algebraic and Topological 
Theories \yr 1985
\pages 391--404 \endref 

\endRefs
\enddocument